\documentclass[11pt]{amsart}
\usepackage{amsmath,amsfonts,amsthm,mathrsfs,amssymb,cancel}
\usepackage[dvipsnames,table,xcdraw]{xcolor}
\usepackage{times}
\usepackage{relsize}
\usepackage[foot]{amsaddr}
\usepackage{tikz}
\usetikzlibrary{cd,arrows}
	\tikzset{thick/.style={line width=.4mm}}
	\tikzstyle{dot}=[circle,thick,fill=red!60]
	\tikzstyle{tinydot}=[circle,thick,fill=red!60,inner sep=1mm]
\usepackage[linktocpage=true]{hyperref}
\hypersetup{
    hypertexnames=true,
    colorlinks=true,
    linkcolor=blue,  
    urlcolor=blue,
    citecolor=blue
}

\newcommand{\colim}{\operatorname{colim}}

\theoremstyle{definition}
\newtheorem{definition}{Definition}
\newtheorem*{definition*}{Definition}

\theoremstyle{theorem}
\newtheorem{lemma}{Lemma}
\newtheorem{theorem}{Theorem}
\newtheorem{corollary}{Corollary}

\title{An enriched category theory of language:  from syntax to semantics}
\author{Tai-Danae Bradley$^1$}
\address{$^1$X, The Moonshot Factory and Sandbox@Alphabet, Mountain View, CA}
\email{tai.danae@math3ma.com}
\author{John Terilla$^2$}
\address{$^2$The City University of New York and Tunnel, New York, NY}
\email{jterilla@gc.cuny.edu}
\author{Yiannis Vlassopoulos$^3$}
\address{$^3$Tunnel, New York, NY}
\email{yiannis@tunnel.tech}

\begin{document}

\begin{abstract}
State of the art language models return a natural language text continuation from any piece of input text.  This ability to generate coherent text extensions implies significant sophistication, including a knowledge of grammar and semantics.  In this paper, we propose a mathematical framework for passing from probability distributions on extensions of given texts, such as the ones learned by today's large language models, to an enriched category containing semantic information.  Roughly speaking, we model probability distributions on texts as a category enriched over the unit interval. Objects of this category are expressions in language, and hom objects are conditional probabilities that one expression is an extension of another.  This category is syntactical---it describes what goes with what. Then, via the Yoneda embedding, we pass to the enriched category of unit interval-valued copresheaves on this syntactical category.  This category of enriched copresheaves is semantic---it is where we find meaning, logical operations such as entailment, and the building blocks for more elaborate semantic concepts.
\end{abstract}
\maketitle
\tableofcontents

\section{Introduction}
The world's best large language models (LLMs) have recently attained new levels of sophistication by effectively learning a probability distribution on possible continuations of a given text.  Interactively, one can input prefix text and then sample repeatedly from a next word distribution to generate original, high quality texts \cite{vaswani2017attention, GPT1, GPT2, GPT3}.  Intuitively, the ability to continue a story implies a great deal of sophistication.  A grammatically correct continuation requires a mastery of syntax, careful pronoun matching, part of speech awareness, a sense of tense, and much more.  A language model that effectively learns a probability distribution on possible continuations must apparently also have learned some \emph{semantic} knowledge.  For the continuation of a story to be reasonable and internally consistent requires knowledge of the world:  dogs are animals that bark, golf is played outdoors during the day, Tuesday is the day after Monday, etc. What is striking is that these LLMs can be trained using unlabeled samples of text to predict a next word.  No grammatical or semantic input is provided, nevertheless complex syntactic structures, semantic information, and world knowledge are learned and demonstrated.  The present work is a response to the real-world evidence that it is possible to pass from probability distributions on text continuations to semantic information.  We propose a mathematical framework for this passage.

We define a syntax category, which is a category enriched over the unit-interval $[0,1]$, that models probability distributions on text continuations \cite{BV2020}.  Our semantic category is then defined to be the enriched category of $[0,1]$-valued copresheaves on the syntax category.  The Yoneda embedding, which maps the syntax category as a subcategory of the semantic category, assigns to a given text its representable copresheaf.  We regard the $[0,1]$-valued copresheaf represented by a text as the meaning of the text, as in dynamic semantics \cite{sep-dynamic-semantics}, whose slogan is ``meaning is context change potential.'' 
Furthermore, there are categorical operations in the semantic category that allow one to combine meanings that correspond to certain logical operations.  In particular, there is a kind of context-sensitive implication that models whether a certain text is true, given that another text is true.  

The paper is organized as follows. Section \ref{ssec:composionality} motivates the category theoretical approach to language by contrasting it to an algebro-geometric perspective and describing a few advantages with explicit examples in Section \ref{ssec:why_ct}. As noted there, a primary advantage is the ability to marry both the compositional and distributional structures of language in a principled way, a claim fully developed in Section \ref{sec:enriched_cat}. There, a few basic definitions from enriched category theory are recalled before defining the language syntax category $\mathcal{L}$ as a category enriched over $[0,1]$.  We then review more enriched category theory and prove the results we need in the special case that the enriching category is the unit interval.  This sets the stage to pass to our semantic category of $[0,1]$-valued copresheaves on $\mathcal{L}$.  The next section details operations on enriched copresheaves that are akin to conjunction, disjunction, and implication between meanings of expressions in language.  Section \ref{sec:metric_space} makes use of an isomorphism between $[0,1]$ and the set of nonnegative extended reals $[0,\infty]$ to recast our ideas in the language of generalized metric spaces and tropical geometry. Finally, a conclusion and summary is provided in Section \ref{sec:conclusion}.

The intended audience for this work includes mathematicians as well as mathematical data scientists who may be intrigued by recent advances in deep learning and natural language, as well as researchers interested in interpretability.  Some of the content may be helpful to philosophers concerned with reasoning, artificial intelligence, and logic in both abstract theory and applications.  We do assume the reader is familiar with basic concepts in category theory: categories, functors, natural transformations, limits, and colimits. Excellent introductions to the subject are readily available, including \cite{riehl2017category,leinster2014basic,fong2019invitation,lawvere2009conceptual}. We do not assume familiarity with enriched category theory and introduce definitions and prove results as needed. For a gentle introduction to the topic, see \cite[Chapter 2]{fong2019invitation} as well as \cite{kelly1982basic} for a more formal treatment. For a historical background of statistical language models and common techniques, see the survey in \cite{jing2019}, and for nontechnical discussions of the current best-in-class LLMs see \cite{NYT2020,MIT2020,MIT2021}.

\subsection{Compositionality}\label{ssec:composionality}
Language is a compositional structure:  expressions can be combined to make longer expressions.  This is not a new idea.  Theoretical linguists have studied grammatical structures for a long time --- think, for example, of the richly developed field of formal grammars \cite{enwiki:1032239019}.  Focusing on the compositionality, one may consider language as some kind of algebraic structure, perhaps simply as the free monoid on some set of atomic symbols quotiented by the ideal of all expressions that are not grammatically correct; or, perhaps, as something more operadic, resembling the tree-like structures of parsed texts, labelled with classes that can be assembled by some set of tree-grafting rules.  In either case, the algebraist might identify the meaning of a word, such as ``red'' with the principal ideal of all expressions that contain ``red.''  The concept of red, then, is the ideal containing \emph{red ruby, bright red ruby, red hot chili peppers, red rover, blood red, red blooded, the workers' and peasants' red army, red meat, red firetruck,...}  If one knows every expression that contains a word, then, as the thinking goes, one understands the meaning of that word.  One is then led to a geometric picture by an algebro-geometric analogy: construct a space whose points are ideals and view the algebraic structure as a sort of coordinate algebra on the space.  From this perspective, language is a coordinate algebra on a space of meanings.  This geometric picture points in appealing directions.  If one can say the same kinds of things in different languages, then those languages should be regarded as having comparable ideal structures (or possibly equivalent categories of modules), translating between different languages involves a change of coordinates on the space of meanings, and so on.  Going further, one has a way to model semantics via sheaves on this space of meanings, as Lawvere in his classic 1969 paper \cite{lawvere69} about syntax-semantic duality. 

We think this algebro-geometric picture is inspirational, but incomplete.  Here, we furnish it with an important refinement.  The compositional structure of language is only part of the structure that exists in language.  What is missing is the \emph{distribution of meaningful texts}.  What has been, and what will be, said and written is crucially important.  When one encounters the word ``firetruck'' it's relevant that ``red firetruck'' has been observed more often than ``green firetruck.''  The fact that ``red firetruck'' is observed more often than ``red idea'' contributes to the meaning of \emph{red}.  Encountering an unlikely expression like ``red idea'' involves a shift in expectation that demands attention and contributes to the meaning of the broader context containing it.  Again, the idea that distributional structure is important in language is not new.  The so called \emph{distributional hypothesis} \cite{harris54} in linguistic states that linguistic items with similar distributions have similar meanings.  In machine learning, word frequency counts have been employed to automatically extract knowledge from a text corpus in Latent Semantic Analysis and vector models of meaning \cite{Turney_2010}.

So, the ideals mentioned above are just a first approximation for the space of meanings.  In this paper, we formalize a marriage of the compositional structure of language with the distributional structure as a category enriched over the unit interval, which we denote by $\mathcal{L}$ and call \emph{the language syntax category}.  Then, once we have defined this language syntax category, we pass to the category of enriched copresheaves on $\mathcal{L}$, the objects of which are unit interval valued functions on $\mathcal{L}$ satisfying a certain monotonicity condition.  The category $\widehat{\mathcal{L}}$ of enriched copresheaves, which is itself a category enriched over the unit interval, is the \emph{semantic} category.  Meanings reside in this semantic category, which also admits ways for meanings to be manipulated and combined to form higher semantic concepts.

Before we discuss why we are using category theory, let us comment briefly on two other efforts to study natural language using categorical methods.
The DisCoCat program of \cite{clark2010} seeks to combine compositional and distributional structures in a single categorical framework, though a choice of grammar is needed as input.  In our work, syntactical structures are inherent within the enriched category theory and no grammatical input is required. This is motivated by the observation that LLMs can continue texts in a grammatically correct way without the additional input of a grammatical structure, suggesting that all important grammatical information is already contained in the enriched category.  In short, the DisCoCat program aims to attribute meaning to parts of texts and grammatical rules for combining them to build meaning of larger texts.  This is like the reverse of our work, which asserts that the meaning of small texts is derived from the distribution of larger texts that contain them.  In \cite{Abramsky2014}, Abramsky and Sadrzadeh consider a sheaf theoretic framework for studying language.  That work makes some interesting use of the gluing condition for sheaves.  Here, we work with copresheaves, without any kind of gluing conditions, but the most important difference is that Abramsky and Sadrzadeh are not using the distributional structure that we are focused on here.

\subsection{Why category theory?}\label{ssec:why_ct}
The algebraic perspective of viewing ideals as a proxy for meaning is consistent with a category theoretical perspective, and the latter provides a better setting in which to merge the compositional and distributional structures of language.  But even before adding distributional structures, moving from an algebraic to a categorical perspective provides certain conceptual advantages that we highlight first.

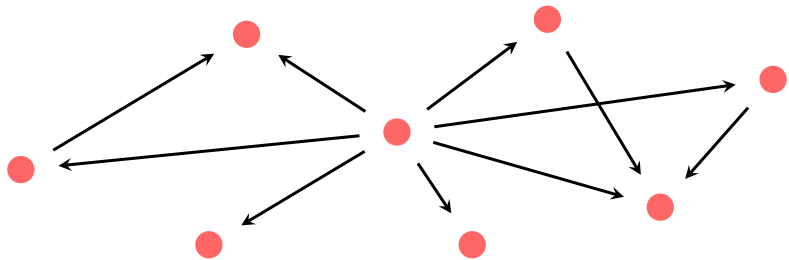
\begin{figure}
    \centering
    \begin{tikzpicture}
	\node[dot] (R) at (0,0) {};
	\node[dot] (A) at (-2,1.3) {};
	\node[dot] (B) at (5,0.7) {};
	\node[dot] (C) at (2,1.5) {};
	\node[dot] (D) at (-2.5,-1.5) {};
	\node[dot] (E) at (-5,-0.5) {};
	\node[dot] (F) at (3.5,-1) {};
	\node[dot] (G) at (1,-1.5) {};

	\draw[thick,-stealth,shorten <=3mm,shorten >=3mm] (R) -- (A) node[midway,above] {};
	\draw[thick,-stealth,shorten <=3mm,shorten >=3mm] (R) -- (B) node[midway,above] {};
	\draw[thick,-stealth,shorten <=3mm,shorten >=3mm] (R) -- (C) node[midway,above] {};
	\draw[thick,-stealth,shorten <=3mm,shorten >=3mm] (C) -- (F) node[midway,above] {};
	\draw[thick,-stealth,shorten <=3mm,shorten >=3mm] (R) -- (D) node[midway,above] {};
	\draw[thick,-stealth,shorten <=3mm,shorten >=3mm] (R) -- (E) node[midway,above] {};
	\draw[thick,-stealth,shorten <=3mm,shorten >=3mm] (R) -- (F) node[midway,above] {};
	\draw[thick,-stealth,shorten <=3mm,shorten >=3mm] (R) -- (G) node[midway,above] {};
	\draw[thick,-stealth,shorten <=3mm,shorten >=3mm] (B) -- (F) node[midway,above] {};
	\draw[thick,-stealth,shorten <=3mm,shorten >=3mm] (E) -- (A) node[midway,above] {};
    \end{tikzpicture}
    \caption{A category with at most one morphism between any two objects. Here the objects are expressions in a language and morphisms indicate when one expression is a continuation of another.  Idenity morphisms are not pictured, but are understood to be present.}
    \label{fig:dots}
\end{figure}

Consider the category whose objects are elements of the free monoid on some finite set of atomic symbols, where there is a morphism $x\to y$ whenever $x$ is a substring of $y$, that is, when the expression $y$ is a \emph{continuation} of $x$. If the finite set is taken to be a set of English words, for example, then there are morphisms \emph{red} $\to$ \emph{red firetruck} and \emph{red} $\to$ \emph{bright red ruby} and so on. Each string is a substring of itself,  providing the identity morphisms, and composite morphisms are provided by transitivity: if $x$ is a substring of $y$ and $y$ is a substring of $z$, then $x$ is a substring of $z.$ This category is simple to visualize: it is \emph{thin}, which is to say there is at most one morphism between any two objects, and so one might have in mind a picture like that in Figure \ref{fig:dots}.  A consequence of the Yoneda lemma implies that a fixed object in this category is  determined up to unique isomorphism by the totality of its relationships to all other objects in the category. One thus thinks of identifying an expression $x$ with the functor $h^x:=\hom(x,-)$ whose value on an expression $y$ is the one-point set $\ast$ if $x$ is contained in $y$ and is the empty set otherwise. The functor $h^x$ is an example of a \emph{copresheaf}, the name given to a functor from a given category to the category of sets, and is in fact a \emph{representable copresheaf}, represented by the object $x$.  So in this way the preimage of $\ast$ is precisely the principal ideal generated by $x.$  Representable copresheaves, therefore, are comparable to principal ideals and both can be thought of as a first approximation for the meaning of an expression.

\subsection{Constructions in the unenriched setting} 
An immediate advantage to the shift in perspective is that the functor category of copresheaves $\widehat{\mathsf{C}}:=\mathsf{Set}^{\mathsf{C}}$ on any small category $\mathsf{C}$ is better behaved than $\mathsf{C}$ itself. It has all limits and colimits and is Cartesian closed.  In fact, the category of copresheaves is a \emph{topos}, which is known to be the appropriate setting for intuitionistic logic \cite{maclane2012sheaves}.  
The full theory of toposes is not incorporated in this work, but intuition abounds when one takes inspiration from the field.  Importantly, one has concrete operations for building new copreasheaves from representable ones, suggestive of the ``meanings are composable instructions'' perspective of internalist semantics \cite{Pietroski}.  

\subsubsection{Coproducts}\label{sssec:coprod}  
In the algebraic setting, the union of ideals is not an ideal, and so one is left to wonder what algebraic structure might represent the disjunction of two concepts, say ``red'' or ``blue.'' In the categorical setting, the answer is straightforward. The coproduct of copresheaves is again a copreasheaf. Given any objects $x$ and $y$ in a small category $\mathsf{C}$, the coproduct of representable copresheaves $h^x \sqcup h^y\colon \mathsf{C}\to\mathsf{Set}$ is computed ``pointwise.'' So suppose $\mathsf{C}$ is the category of expressions in language defined above, and let $x=\textit{red}$ and $y=\textit{blue}$. The coproduct $h^\textit{red} \sqcup h^\textit{blue}$ therefore maps an expression $c$ to the set $h^\textit{red}(c) \sqcup h^\textit{blue}(c)$. This set is isomorphic to $\ast$ if $c$ contains either \emph{red} or \emph{blue}, and it is isomorphic to a two-point set if $c$ contains both, and otherwise it is the empty set. The output of the functor $h^\textit{red}\sqcup h^\textit{blue}$ is therefore nonempty on the union of all expressions that contain \emph{red} or \emph{blue} or both, and this matches well with the role of union as logical ``or'' among sets.

\subsubsection{Products}\label{sssec:prod}
Now think about limits and in particular products. Like the coproduct, the product of representable copresheaves $h^x\times h^y\colon\mathsf{C}\to\mathsf{Set}$ is again a copresheaf computed pointwise. So if $\mathsf{C}$ is the category of language and $x=\textit{red}$ and $y=\textit{blue}$, then the value of the product on any expression $c$ is given by $h^\textit{red}(c)\times h^\text{blue}(c)$, which is isomorphic to $\ast$ if $c$ contains both \emph{red} and \emph{blue} and is the empty set otherwise. So the output of the functor  $h^\textit{red}\times h^\textit{blue}$ is nonempty on the intersection of expressions that contain \emph{red} with those that contain \emph{blue}, and this coincides with the role of intersection as logical ``and'' among sets.

\subsubsection{Cartesian closure}\label{sssec:implication}
Advantages of the categorical perspective can further be illustrated with a third example, which comes from the Cartesian closure of the category of copresheaves $\widehat{\mathsf{C}}$ on a small category $\mathsf{C}$. To say that $\widehat{\mathsf{C}}$ is \emph{Cartesian closed} means that it has finite products and moreover that the product of copresheaves has a right adjoint (See section I.6.~of \cite{MacLane:824105}).  More precisely, there is a bifunctor $[\,\,,\,\,]\colon \widehat{\mathsf{C}}\times \widehat{\mathsf{C}} \to \widehat{\mathsf{C}}$ called the \emph{internal hom} that fits into an isomorphism $\widehat{\mathsf{C}}(F\times G,H)\cong \widehat{\mathsf{C}}(F,[G,H])$ for all copresheaves $F,G$ and $H$. Here and throughout we use the notation $\mathsf{A}(x,y)$ to denote the set of morphisms from objects $x$ to $y$ in a category $\mathsf{A}$.
The relevance of this new copresheaf $[G,H]$ can be seen in its kinship to \emph{implication} in logic. Indeed, in intuitionistic propositional calculus one works within a certain kind of thin Cartesian closed category called a \emph{Heyting algebra}.  The usual notation in this setting is to write $x\leq y$ if there is a morphism from $x\to y$.  The categorical product is denoted by $\wedge$ and the internal hom is denoted by $\Rightarrow$.  So, in this notation, one has: $x\wedge y\leq z$ if and only if $x\leq (y\Rightarrow z)$ for all elements $x,y,z$. Keep this analogy in mind as we dissect an example of the internal hom in $\widehat{\mathsf{C}}$ when $\mathsf{C}$ is the category of expressions in language.  Let's do some unwinding to see that the copresheaf \[[h^\textit{red},h^\textit{blue}] \colon \mathsf{C} \to \mathsf{Set}\] captures something like an implication \emph{red}~$\Rightarrow$~\emph{blue}.  First use the Yoneda lemma and then the defining property of the internal hom to get
\begin{equation}\label{eq:cartclosed}
    [h^\textit{red},h^\textit{blue}](c) = \widehat{\mathsf{C}}(h^c,[h^\textit{red},h^\textit{blue}]) = \widehat{\mathsf{C}}(h^c\times h^\textit{red},h^\textit{blue}).
\end{equation}
So,  the internal hom assigns to an expression $c$ the set of natural transformations from the functor $h^c\times h^\textit{red}$ to the functor $h^\textit{blue}$.  The data of a single natural transformation from $h^c\times h^\textit{red}$ to $h^\textit{blue}$ consists of a collection of set functions
\begin{equation}\label{eq:implicationpresheaf}
\left\{(h^c\times h^{\textit{red}})(d)\to h^\textit{blue}(d)\right\}_{d\in  \mathrm{ob}(\mathsf{C})}
\end{equation} that fit into a particular commutative square.  Here, the domain $(h^c\times h^{\textit{red}})(d)$ and the codomain $(h^{\textit{blue}})(d)$ are either empty or a singleton, and so the data of a natural transformation either does not exist or is uniquely specified and automatically fits into the required commutative square.  Therefore, $[h^\textit{red},h^\textit{blue}](c)$ is either the empty set $\emptyset$ or the one-point set $\ast$. 

To determine whether $[h^\textit{red},h^\textit{blue}](c)$ is empty or not, let's look closer at the functions in \eqref{eq:implicationpresheaf}.  The product is computed pointwise, and thus the domain $(h^c\times h^{\textit{red}})(d)=h^c(d)\times h^{\textit{red}}(d)$ is isomorphic to the one-point set $\ast$ whenever $d$ contains both $c$ and \emph{red} and is the empty set otherwise. The codomain  $h^\textit{blue}(d)$ is likewise either $\ast$ if $d$ contains \emph{blue} and is the empty set otherwise. Note that if there exists a text $d$ that contains $c$ and \emph{red} but does not contain \emph{blue} then $(h^c\times h^{\textit{red}})(d) =\ast$ and $h^\textit{blue}(d)=\emptyset$, hence there does not exist a function$(h^c\times h^{\textit{red}})(d) \to h^\textit{blue}(d)$ and the set of natural transformations $h^c\times h^\textit{red}$ to $h^\textit{blue}$ is empty.  On the other hand, if every text $d$ that contains $c$ and \emph{red} also contains \emph{blue} then there is a unique function $(h^c\times h^{\textit{red}})(d)\to h^\textit{blue}(d)$ which is specified as follows:
\[
\begin{cases}
\ast \to \ast & \text{when $d$ contains $c$, \emph{red}, and \emph{blue}}\\
\emptyset \to \ast & \text{when $d$ does not contain both $c$ and \emph{red}, and $d$ does contain \emph{blue}}\\
\emptyset \to \emptyset & \text{when $d$ does not contain both $c$ and \emph{red}, and $d$ does not contain \emph{blue}}.
\end{cases}
\]
For example, if $c$ is the expression \emph{French flag}, then there exists a unique natural transformation from $h^c \times h^\textit{red}$ to $h^\textit{blue}$ if every text that contains both \emph{French flag} and \emph{red} as subtexts, also contains $\emph{blue}$ and so $[h^\textit{red},h^\textit{blue}](\emph{French flag}) = \ast$.  If $c$ is the expression \emph{ruby} then there exists no natural transformation $h^c \times h^\textit{red}$ to $h^\textit{blue}$ if there is a text that contains \emph{red} and \emph{ruby} and does not contain \emph{blue} and so $[h^\textit{red},h^\textit{blue}](\emph{ruby}) = \emptyset.$  One might say ``\emph{red} implies \emph{blue} in the context of \emph{French flag}'' but ``\emph{red} does not implies \emph{blue} in the context of \emph{ruby}.''


In summary, the copresheaf $[h^x,h^y]\colon \mathsf{C} \to \mathsf{Set}$ is given by
\begin{equation*}
[h^x,h^y](c) = 
\begin{cases}
\ast & \text{ if every text that contains both $c$ and $x$ also contains $y$}\\
\emptyset & \text{otherwise}.
\end{cases}
\end{equation*}


Let's review the picture of language presented in this section. 
One has a simplified language category $\mathsf{C}$ whose objects are expressions in the language having a single morphism $x\to y$ if $x$ is a subexpression of $y$. The functor $\mathsf{C}^{\text{op}}\to\mathsf{Set}^{\mathsf{C}}$ that maps an expression $x$ to the representable functor $h^x=\hom(x,-)$ is called the \emph{Yoneda embedding} and captures something of the meaning of $x$. Here $\mathsf{C}^{\text{op}}$ is the \emph{opposite} category of $\mathsf{C}$. As a category, it has the same objects as $\mathsf{C}$ and the morphisms are defined to be $\mathsf{C}^{\text{op}}(x,y):= \mathsf{C}(y,x).$) This is comparable to making a passage from \emph{syntax} to \emph{semantics}, after which meanings can be combined by computing products, coproducts, and internal homs.  Even more is possible, for combined meanings can be combined again forming higher concepts, and there are other categorical limits and colimits beyond products and coproducts, such as pushouts, pullbacks, equalizers, and so on....  Yet this picture is still incomplete. The distributional structure of language has not yet been accounted for. For this, we use enriched category theory.


\section{Enriched category theory}\label{sec:enriched_cat}
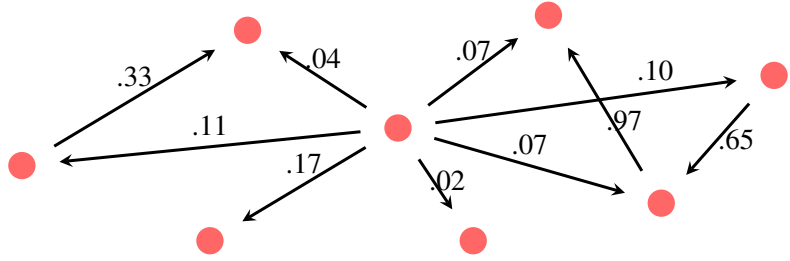
\begin{figure}[h!]
    \centering
    \begin{tikzpicture}
	\node[dot] (R) at (0,0) {};
	\node[dot] (A) at (-2,1.3) {};
	\node[dot] (B) at (5,0.7) {};
	\node[dot] (C) at (2,1.5) {};
	\node[dot] (D) at (-2.5,-1.5) {};
	\node[dot] (E) at (-5,-0.5) {};
	\node[dot] (F) at (3.5,-1) {};
	\node[dot] (G) at (1,-1.5) {};

	\draw[thick,-stealth,shorten <=3mm,shorten >=3mm] (R) -- (A) node[midway,above] {.04};
	\draw[thick,-stealth,shorten <=3mm,shorten >=3mm] (R) -- (B) node[pos=0.7,above] {.10};
	\draw[thick,-stealth,shorten <=3mm,shorten >=3mm] (R) -- (C) node[midway,above] {.07};
	\draw[thick,-stealth,shorten <=3mm,shorten >=3mm] (F) -- (C) node[pos=0.3,above] {.97};
	\draw[thick,-stealth,shorten <=3mm,shorten >=3mm] (R) -- (D) node[midway,above] {.17};
	\draw[thick,-stealth,shorten <=3mm,shorten >=3mm] (R) -- (E) node[midway,above] {.11};
	\draw[thick,-stealth,shorten <=3mm,shorten >=3mm] (R) -- (F) node[midway,above] {.07};
	\draw[thick,-stealth,shorten <=3mm,shorten >=3mm] (R) -- (G) node[pos=0.7,above] {.02};
	\draw[thick,-stealth,shorten <=3mm,shorten >=3mm] (B) -- (F) node[pos=0.3,below] {.65};
	\draw[thick,-stealth,shorten <=3mm,shorten >=3mm] (E) -- (A) node[midway,above] {.33};
    \end{tikzpicture}
    \caption{The compositional and distributional structures of language are married by decorating arrows with the conditional probability that one expression contains another.}
    \label{fig:enriched_dots}
\end{figure}

Enriched category theory provides a ready-made way to decorate morphisms in the simplified language category $\mathsf{C}$ from Section \ref{ssec:why_ct} with conditional probabilities as in Figure \ref{fig:enriched_dots}.  Enriched category theory begins with the observation that the set of morphisms between objects may have more structure than just a \emph{set}.  Examples are plentiful. The set of linear maps between vector spaces is itself a vector space, for instance.  In enriched category theory, the morphisms between objects is itself an object in a category called the \emph{base category} or the category over which the category is \emph{enriched}.  When the enriching category is the category $\mathsf{Set}$ of sets, enriched category theory reduces to ordinary category theory.  In order for enriched category theory to have the desired structures and axioms (such as having composable morphisms), the base category must have some of the structures that the category of sets has.  One can go rather far assuming the base category is a symmetric monoidal category.  In order to have convenient versions of enriched presheaves and copresheaves, the base category should be a symmetric monoidal category that is also \emph{closed}, which allows the base category to be enriched over itself.  If, in addition, the base category is \emph{complete} and \emph{cocomplete}, meaning that it contains all limits and colimits, then the categories of copresheaves and presheaves are complete and cocomplete.  From a categorical point of view, our setting is relatively simple.  The  set $\mathsf{C}(x,y)$ of morphisms from $x$ to $y$ is either the empty set or the one point set.  Also, the category we wish to enrich over is the unit interval $[0,1]$, which becomes a complete, closed, symmetric monoidal category in a simple way that fits our purposes well.  So, we will not burden the reader with the full machinery of general enriched category theory \cite{kelly1982basic} but rather take advantage of the simplifications afforded by our setting and specialize some of the definitions.

\subsection{Categories enriched over $[0,1]$}\label{sec:enrichedlangcat}

A \emph{preorder} is a set together with a reflexive, transitive relation.  Any preorder becomes a category whose objects are the elements of the set, with a single morphism from $a$ to $b$ if and only if $a$ is related to $b$.  Reflexivity provides identity morphisms and transitivity provides composition of morphisms, which is defined in the only way it can be.  Limits and colimits of finite diagrams always exists and are easy to compute:  the limit is the minimum of the elements in the diagram and the colimit of a diagram is the maximum.

A \emph{monoid} is a set together with an associative binary operation and a unit for that operation.   A \emph{commutative} monoid is a monoid whose operation is commutative.  A preorder with a compatible commutative monoidal structure naturally becomes a kind of category over which one can enrich.  Formally,

\begin{definition}
A \emph{commutative monoidal preorder} $(\mathcal{V},\leq,\otimes, 1)$ is a preorder $(\mathcal{V}, \leq)$ and a commutative monoid $(\mathcal{V},\otimes,1)$ satisfying $x\otimes y \leq x'\otimes y'$ whenever $x\leq x'$ and $y \leq y'$.  
\end{definition}

\begin{definition}
Let $(\mathcal{V},\leq,\otimes,1)$ be a commutative monoidal preorder.  The data of a (small) \emph{$\mathcal{V}$-enriched category}, or simply a \emph{$\mathcal{V}$-category}, consists of a set of objects $\mathcal{C}$, and for every pair of objects $x$ and $y$, there is an element $\mathcal{C}(x,y)\in\mathcal{V}$ called a \emph{$\mathcal{V}$-hom object}.  This data satisfies  
\begin{gather}\label{enrichedcateq1}
1\leq \mathcal{C}(x,x) \\
\mathcal{C}(y,z)\otimes\mathcal{C}(x,y)\leq\mathcal{C}(x,z)\label{enrichedcateq2}
\end{gather} 
for all objects $x,y,z \in \mathcal{C}$.
\end{definition}

The unit interval $[0,1]:=\{x\in \mathbb{R}:0\leq x \leq 1\}$ is a commutative monoidal preorder with multiplication being the monoidal product, having $1$ as the unit, and with the usual $\leq$ relation being the preorder.  The data of a $[0,1]$-category consists of a set of objects $\mathcal{C}$ and a  $[0,1]$-valued function $(x,y) \mapsto \mathcal{C}(x,y)$ defined for every $x,y \in \mathcal{C}$ satisfying 
$\mathcal{C}(x,x)=1$ for every $x\in \mathcal{C}$
and $\mathcal{C}(y,z)\mathcal{C}(x,y) \leq \mathcal{C}(x,z)$ for all $x,y,z \in \mathcal{C}$.

\begin{definition}
A commutative monoidal preorder $\mathcal{V}$ is said to be \emph{closed} provided for every pair of elements $x$ and $y$ in $\mathcal{V}$ there is an element $[x,y]\in \mathcal{V}$ , called the \emph{internal hom}, satisfying 
\begin{equation}\label{eq:closed}
x\otimes y \leq z \text{ if and only if }x \leq [y,z]
\end{equation} for all $x,y,z$ in $\mathcal{V}$. 
\end{definition}

The relevance is that a \emph{closed} commutative monoidal preorder $\mathcal{V}$ becomes a category enriched over itself by 
replacing the hom set $\mathcal{V}(x,y)$, which is either the empty set or the one-point set, by the internal hom $[x,y]$, which is an object of $\mathcal{V}$.  To see that the assignment $(x,y)\mapsto [x,y]$ does in fact make $\mathcal{V}$ a category enriched over itself, one needs to check that $1\leq [x,x]$ and $[y,z]\otimes[x,y] \leq [x,z]$.  The first inequality follows from the fact that $1\otimes x \leq x$ and Equation \eqref{eq:closed}.  One can check that $[y,z]\otimes [x,y] \leq [x,z]$ in two steps.  First, $[y,z] \leq [y,z]$ and Equation \eqref{eq:closed} gives 
$[y,z]\otimes y \leq z$ and similarly, $[x,y]\otimes x \leq y$.  Putting these two inequalities together yields $[y,z] \otimes [x,y] \otimes x \leq z$ which implies $[y,z] \otimes [x,y]\leq [x,z]$. Much can be said about commutative monoidal preorders and enrichment, for instance, see \cite[Chapter 2]{fong2019invitation}, though for now our primary focus will be on the unit interval.

\begin{lemma}
\label{def:interval}
The unit interval $[0,1]$ is a closed commutative monoidal preorder.  The monoidal product $a\otimes b:=ab$ is the usual product of numbers and the internal hom is \emph{truncated division}: for all $a,b\in[0,1]$ define
\begin{equation}\label{eq:truncated_div}
		[a,b]:=
		\begin{cases}
		b/a &\text{if } b< a,\\
		1 &\text{otherwise}.
		\end{cases}
\end{equation}
\end{lemma}
\begin{proof}
To verify closure, we need to check the formula for truncated division works as an internal hom; that is, $ab \leq c $ if and only if $a \leq [b,c]$.  There are two cases.  If $c<b$ then $ab\leq c$ implies that $a \leq \frac{c}{b}=[b,c].$  If $b\leq c$ then $[b,c]=1$ and $a \leq [b,c]$ automatically.
\end{proof}
While everything constructed using the internal hom in $[0,1]$ ultimately involves statements about multiplication and less than or equal to, trying to argue using only multiplication and order can get messy, as statements often break down into multiple cases owing to the minimum in truncated division.  Using the fact that for all $a,b,c \in [0,1]$ we have
\begin{equation}\label{eq:Iclosed}
ab \leq c \text{ if and only if }a \leq [b,c]
\end{equation} 
often makes arguments simpler.  In what follows, we refer to the equivalence in \eqref{eq:Iclosed} as the \emph{closure} property of $[0,1]$.  We occasionally will write $\min\{b/a,1\}$ when $a$ might be zero.  The reader should interpret $b/0 \geq 1$ for any $0\leq b \leq 1.$

Also, keep in mind that we use juxtaposition $ab$ for ordinary multiplication of numbers.  For our purposes, this is shorthand for the \emph{monoidal} product $a\otimes b.$  We also have the \emph{categorical product} in $[0,1]$, which is given by the minimum.  So, for two numbers $a,b\in [0,1]$ the expression $a\times b$ means $\min\{a,b\}.$ The \emph{coproduct} is given by the maximum and is denoted by $a\sqcup b$.  So, when we are working in the unit interval $[0,1]$, remember Equation \eqref{eq:truncated_div} for the internal hom and 
\begin{align*}
a\times b&:= \min\{a,b\} \\
a\sqcup b&:= \max\{a,b\}.
\end{align*}

As a final remark, any thin category $\mathsf{C}$ can become a category enriched over $[0,1]$ by setting $\mathsf{C}(x,y)=1$ if there is a morphism between objects $x\to y$ and $\mathsf{C}(x,y)=0$ otherwise.  In fact, $2:=\{0,1\}$, the set that contains zero and one only, is itself a commutative monoidal preorder that is closed, a sub-commutative monoidal preorder of the unit interval, and can serve as a base category over which to enrich.

\subsection{The syntax category $\mathcal{L}$}

Following \cite{BV2020}, we now define a category $\mathcal{L}$ enriched over $[0,1]$.
\begin{definition}\label{def:syntax}We define the \emph{syntax category} $\mathcal{L}$ to be the category enriched over $[0,1]$ whose objects are expression in the language and where $[0,1]$-objects are defined by 
\[\mathcal{L}(x,y):=\pi(y\vert x),\] 
where $\pi(y\vert x)$ denotes the probability that expression $y$ extends expression $x$. If $x$ is not a subtext of $y$ then necessarily $\pi(y\vert x)=0$.  
\end{definition}
So, for example, one might have 
\begin{align*}
\mathcal{L}\left(\textit{red,red firetruck}\right) &=0.02\\
\mathcal{L}\left(\textit{red,red idea}\right) &=10^{-5}\\
\mathcal{L}\left(\textit{red,blue sky}\right) &=0.
\end{align*}
That $\mathcal{L}$ is indeed a category enriched over $[0,1]$ follows from the fact that $\pi(x|x)=1$ and
\begin{equation}\label{eq:Lenriched}
\pi(z \vert y)\pi(y\vert x)=\pi(z\vert x)
\end{equation}
for all texts $x,y,z$ and so satisfies the required inequalities \eqref{enrichedcateq1} and \eqref{enrichedcateq2} with equalities.  The reader might think of these probabilities $\pi(y \vert x)$ as being most well defined when $y$ is a short extension of $x$.  While one may be skeptical about assigning a probability distribution on the set of all possible texts, it's reasonable to say there is a nonzero probability that \emph{cat food} will follow \emph{I am going to the store to buy a can of} and, practically speaking, that probability can be estimated.  Indeed, existing LLMs \cite{GPT1, GPT2, GPT3} successfully learn these conditional probabilities $\pi(y\vert x)$ using standard machine learning tools trained on large corpora of texts, which may be viewed as providing a wealth of samples drawn from these conditional probability distributions.  Figure \ref{fig:enriched_dots} gives the right toy picture:  the objects are expressions in language, and the labels on the arrows describe the probability of extension.   

As in the unenriched setting of Section \ref{ssec:why_ct}, the category $\mathcal{L}$ is inherently \emph{syntactical}---it encodes ``what goes with what'' together with the statistics of those expressions.  But what about \emph{semantics}? Following the same line of reasoning that led to the consideration of copresheaves in the unenriched case, we now wish to pass from $\mathcal{L}$ to the enriched version of \emph{copresheaves on $\mathcal{L}$} where we propose that meaning lies, and where concepts can again be combined through the enriched version of limits, colimits, and internal homs. This discussion first requires a few mathematical preliminaries, beginning with what an enriched functor is, what an enriched copresheaf is, and what the enriched version of natural transformations between functors are.

Before going on to the next section, we briefly comment on the work in \cite{BV2020}.  In that paper, a $[0,1]$-enriched functor is defined between the syntax category $\mathcal{L}$ and another $[0,1]$-enriched category, embedding $\mathcal{L}$ as subcategory of a $[0,1]$-enriched category $\mathcal{D}$ consisting of ``density'' operators.  In this paper we construct a category $\widehat{\mathcal{L}}$ of $[0,1]$-copresheaves on $\mathcal{L}$, and the enriched version of the Yoneda embedding defines an embedding ${\mathcal{L}}\to \widehat{\mathcal{L}}$.  The reader may think of the $[0,1]$-enriched Yoneda embedding as factoring through the functor $\mathcal{L}\to\mathcal{D}$ defined in \cite{BV2020}, as in the picture below:
\begin{equation*}
    \begin{tikzcd}
    \mathcal{L} \arrow[rr] \arrow[dr] && \widehat{\mathcal{L}}\\
    & \mathcal{D} \arrow[ru,dotted]\\
    \end{tikzcd}
\end{equation*}
 While technical details are needed to construct the third arrow $\begin{tikzcd}\mathcal{D} \arrow[r,dotted] &\widehat{\mathcal{L}}\end{tikzcd}$, one may think of the  category $\mathcal{D}$ is a kind of intermediary.  The logical and semantic possibilities within $\mathcal{D}$ are more limited than in the semantic $\widehat{\mathcal{L}}$.  However, the category $\mathcal{D}$, as is argued in \cite{BV2020}, can be approximated by a computer model, while providing more room for semantic exploration than in the syntactic category $\mathcal{L}$.  

\section{Enriched copresheaves}\label{sec:enriched_copresheaves}
 Recall from Section \ref{ssec:why_ct} that a copresheaf on a category $\mathsf{C}$ is a functor $\mathsf{C}\to\mathsf{Set}$. To understand the enriched version, then, one must first have a notion of functors between enriched categories. 
 \begin{definition}
 Suppose $\mathcal{C}$ and $\mathcal{D}$ are categories enriched over a commutative monoidal preorder $(\mathcal{V},\otimes,\leq,1)$.  An \emph{enriched functor} $\mathcal{C}~\to~\mathcal{D}$ is a function $f\colon \mathcal{C}\to \mathcal{D}$ satisfying
\begin{equation}\label{eq:v_functor}
    \mathcal{C}(x,y)\leq \mathcal{D}(fx,fy)
\end{equation}
for all objects $x$ and $y$ in $\mathcal{C}$.  
\end{definition}

In the case that $\mathcal{V}$ is closed, it is enriched over itself and we can take $\mathcal{D}=\mathcal{V}$ to make sense of enriched copresheaves.  

\begin{definition}Suppose $\mathcal{C}$ is a category enriched over a closed commutative monoidal preorder $\mathcal{V}$.  An  \emph{enriched copresheaf} is a function $f\colon \mathcal{C}\to \mathcal{V}$ satisfying $\mathcal{C}(x,y)\leq \mathcal{V}(fx,fy)$ for all objects $x$ and $y$ in $\mathcal{C}.$
\end{definition}

Now, we discuss how to make a $\mathcal{V}$-enriched category $\mathcal{D}^{\mathcal{C}}$  whose objects are $\mathcal{V}$-functors from $\mathcal{C}$ to $\mathcal{D}$. We need a hom object $\mathcal{D}^{\mathcal{C}}(f,g)\in \mathcal{V}$ between any two such functors $f,g\colon\mathcal{C}\to\mathcal{D}$.  Formally, it is defined by an \emph{end}, which is a particular kind of limit in $\mathcal{V}$
\begin{equation}\label{eq:int_v_hom}
    \mathcal{D}^{\mathcal{C}}(f,g):=\int_{c\in\mathcal{C}}\mathcal{D}(fc,gc)
\end{equation}
which always exists if $\mathcal{V}$ is complete.  The precise definition of an end can be found in Section 7.3 in \cite{riehl_2014}, but now we will specialize to the case of copresheaves and the right hand side of Equation \eqref{eq:int_v_hom} reduces to something simple, summarized in Lemma \ref{lem:end_computed} below. The takeaway is that just as mappings between (ordinary) categories define a functor category, so mappings between enriched categories define an \emph{enriched} functor category, at least when the base category is closed and complete.  

From now on, we specialize to the case that the enriching category is $[0,1]$.  We begin with a lemma that says how to assign an element of $[0,1]$ to a pair of functors enriched over $[0,1]$.

\begin{lemma}\label{lem:end_computed}
If $\mathcal{C}$ is a category enriched over $[0,1]$, then the category $\widehat{\mathcal{C}}:=[0,1]^\mathcal{C}$ of copresheaves is also enriched over $[0,1]$.  The $[0,1]$-object between any pair of copresheaves $f,g\colon \mathcal{C}\to[0,1]$ is given by the following infimum over all objects in $\mathcal{C}$,
\begin{equation}\label{eq:hom_min}
\widehat{\mathcal{C}}(f,g)=\inf_{c\in\mathcal{C}}\{[fc,gc]\}
\end{equation}
\end{lemma}
\begin{proof}
The proof is a computation of the end on the right hand side of Equation (\ref{eq:int_v_hom}) when the target category is the unit interval.   
\end{proof}
Recalling that the internal hom in the unit interval is given by truncated division, one has the $[0,1]$-hom object associated to a pair of $[0,1]$-functors $f,g\colon \mathcal{C} \to [0,1]$ is given by $\widehat{\mathcal{C}}(f,g)=\inf_{c\in\mathcal{C}}\{1,gc/fc\}.$

\begin{lemma}\label{lem:representables}
Let $\mathcal{C}$ be a $[0,1]$-category.  For every object $x$, the function 
\[h^x:=\mathcal{C}(x,-)\] is a $[0,1]$-functor.
\end{lemma}
\begin{proof}
Setting $\mathcal{V}=[0,1]$ and also $\mathcal{D}=[0,1]$, Equation \eqref{eq:v_functor} says what is required for a function $\mathcal{C} \to [0,1]$ to be an enriched functor.  So, let $c,d \in \mathcal{C}$.  We need to check that 
\begin{equation*}
\label{eq:lemma1inq}\mathcal{C}(c,d)\leq [h^x(c),h^x(d)]
\end{equation*} 
That is, we need to check that $\mathcal{C}(c,d)\leq  [\mathcal{C}(x,c),\mathcal{C}(x,d)]$, which by the closure condition in \eqref{eq:Iclosed} is equivalent to $\mathcal{C}(c,d)\mathcal{C}(x,c)\leq \mathcal{C}(x,d)$, which is satisfied since $\mathcal{C}$ is enriched over $[0,1]$.  
\end{proof}

\begin{definition}For any object $x$ in a $[0,1]$-category $\mathcal{C}$, we call the functor $h^x:=\mathcal{C}(x,-)$ \emph{the copresheaf represented by $x$.}  We say a copresheaf $f\colon\mathcal{C} \to [0,1]$ is \emph{representable} if $f=h^x$ for some $x\in \mathcal{C}$.
\end{definition}

\subsection{The $[0,1]$-enriched Yoneda lemma}
Lemma \ref{lem:end_computed} says that enriched copresheaves form an enriched category.  Lemma \ref{lem:representables} says that each object defines a representable copresheaf.  In this subsection, we see that the assignment $x\mapsto h^x$ defines an enriched functor that embeds (the opposite of) a $[0,1]$-category within its category of copresheaves.  It is a corollary of an enriched Yoneda lemma.

\begin{theorem}[The Enriched Yoneda Lemma]\label{thm:enrichedyoneda}For any object $x$ in a $[0,1]$ category $\mathcal{C}$ and any $[0,1]$-copresheaf $f\colon\mathcal{C}\to[0,1]$ we have $\widehat{\mathcal{C}}(h^x,f)=f(x).$
\end{theorem}
\begin{proof}Fix an object $x\in \mathcal{C}$ and a copresheaf $f$.
Since $\widehat{\mathcal{C}}(h^x,f)=\inf_{c \in \mathcal{C}}\{[h^x(c),fc]\}$ we have for any particular $c\in \mathcal{C}$ that $\widehat{\mathcal{C}}(h^x,f) \leq [h^x(c),fc]$.  For $c=x$, we have \[\widehat{\mathcal{C}}(h^x,f)\leq [h^x(x),fx] = [1,fx] =fx.\]  On the other hand, since $f$ is a $[0,1]$-functor from $\mathcal{C}$ to $[0,1]$, we have $\mathcal{C}(x,c) \leq [fx,fc]$ for any $c\in \mathcal{C}$.  By the closure of $[0,1]$ the inequality $\mathcal{C}(x,c) \leq [fx,fc]$ is equivalent to $\mathcal{C}(x,c)fx \leq fc$ which in turn is equivalent to $fx \leq [\mathcal{C}(x,c),fc]$.  Having $fx \leq [\mathcal{C}(x,c),fc]$ for every $c\in \mathcal{C}$ implies that $fx \leq \inf_{c \in \mathcal{C}}\{[h^x(c),fc]\}=\widehat{\mathcal{C}}(h^x,f)$ and the theorem is proved.
\end{proof}

\begin{corollary}\label{cor:isometric}
$\mathcal{C}(y,x)= \widehat{\mathcal{C}}(h^x,h^y)$ for all objects $x,y$ in a $[0,1]$-category  $\mathcal{C}$.
\end{corollary}
\begin{proof}
Setting $f=h^y$ in Theorem \ref{thm:enrichedyoneda} yields $\widehat{\mathcal{C}}(h^x,h^y) = h^y(x)=\mathcal{C}(y,x).$
\end{proof}
Therefore, we have the expected interpretation of Corollary \ref{cor:isometric} as an enriched version of the (co)Yoneda embedding.
\begin{corollary}\label{cor:enrichedyoneda}For any $[0,1]$-category $\mathcal{C}$, the assignment $x\mapsto h^x$ defines an enriched functor $\mathcal{C}^{\text{op}} \to \widehat{\mathcal{C}}$, embedding $\mathcal{C}^{\text{op}}$ as an enriched subcategory of $\widehat{\mathcal{C}}$.
\end{corollary}
The \emph{op} in $\mathcal{C}^{\text{op}}$ stands for ``opposite'' and is there because the assignment $x\mapsto h^x$ reverses morphisms, as in the statement of Corollary \ref{cor:isometric} above.  

\subsection{The semantic category $\widehat{\mathcal{L}}$}
Often, when a mathematical object $Y$ has nice structure, the set of functions $\{X\to Y\}$ from a fixed $X$ has nice struture also, even if $X$ does not.  Real valued functions on any set form a vector space, etc.... The unit interval has rich structure from the perspective of category theory.  It is commutative monoidal closed and complete and cocomplete.  For any $[0,1]$-enriched category $\mathcal{C}$, the category $\widehat{\mathcal{C}}=[0,1]^\mathcal{C}$ of copresheaves on $\mathcal{C}$, as a category of functors into the unit interval, inherits rich structure, as well.  In particular, it has enriched versions of products, coproducts, and an internal hom.  Corollary \ref{cor:enrichedyoneda}, from this perspective, says that any $[0,1]$-category embeds into the $[0,1]$-category $\widehat{\mathcal{C}}$, which is often much nicer.  

We now look at the copresheaves on the enriched syntax category $\mathcal{L}$ which will provide a place for the combination of concepts in language in a way that's parallel to the ideas explored in Section \ref{ssec:why_ct}.  

\begin{definition}Let $\mathcal{L}$ be the syntax category (Definition \ref{def:syntax}).  The \emph{semantic} category $\widehat{\mathcal{L}}:=[0,1]^\mathcal{L}$ is the $[0,1]$-category of $[0,1]$-enriched copresheaves on the $[0,1]$-category $\mathcal{L}$.  
\end{definition}

For each object $x$ in $\mathcal{L}$, the representable copresheaf $h^x:=\mathcal{L}(x,-)$ is given by the conditional probability of extending $x$,
\begin{equation}\label{eq:enriched_representable}
    c\mapsto h^x(c) :=
    \begin{cases}
    \pi(c\vert x) &\text{if } x\leq c\\
    0 & \text{otherwise}.
    \end{cases}
\end{equation}
where we use ``$x\leq c$'' as shorthand for ``$x$ is contained as a subtext of $c$.'' 

The representable enriched copresheaf $h^x$ is our proposal for the meaning of the expression $x$.  Its support consists of all expressions containing $x$, which coincides with the principal ideal associated to $x$, and $h^x$ further accounts for the statistics associated with those expressions, which is precisely the distributional structure missing from Section \ref{ssec:why_ct}.  In other words, the meaning of a text is the varying potential of all contexts in which it is used, making our definition of $h^x$ as the meaning of the text $x$ consistent with at least several philosophical traditions, including a use theory of meaning \cite{horwich} and dynamic semantics \cite{sep-dynamic-semantics}.   The embedding in Corollary \ref{cor:enrichedyoneda} assigns to a text $x$ its meaning $h^x$.

From a mathematical point of view, the assignment $x\mapsto h^x$ faithfully embeds the syntax category $\mathcal{L}$ into the category of copresheaves on $\mathcal{L}$.  The category $\widehat{\mathcal{L}}$ is complete and cocomplete and so contains all (the enriched versions of) limits and colimits.   The meanings of texts that live in $\widehat{\mathcal{L}}$ can be manipulated and combined to form higher concepts in $\widehat{\mathcal{L}}$, none of which is possible within the confines of $\mathcal{L}$.  These operations on copresheaves are the subject of the next section.

The reversal of arrows in the $[0,1]$-enriched Yoneda Lemma has a nice interpretation here also.  Suppose the text $y$ extends the text $x$ and $\mathcal{L}(x,y)=a \neq 0.$  We have the following picture: 

\[
\begin{tikzcd}
x\arrow[r,"a"] & y &  & {}& h^y \arrow[r,"a"] &h^x
\end{tikzcd}
\] 
On the left, the text $y$ is an extension of the text $x$ in the syntax category $\mathcal{L}$.  Passing to the semantic category $\widehat{\mathcal{L}}$ on the right, $h^x$ is the meaning of the text $x$, which represents the varying potential contexts in which $x$ might appear, and $h^y$ represents the varying potential contexts in which the text $y$ appears.  The contexts that $x$ can appear extend the contexts that in which $y$ can appear.  Continuing an expression restricts the potential contexts in which the expression can be used.


\section{Enriched products and coproducts in $\widehat{\mathcal{L}}$}\label{sec:enriched_prod_coprod}

Think back to Section \ref{sssec:coprod}, for instance, where we described the coproduct of ordinary copresheaves associated to the words \emph{red} and \emph{blue}, which represented their disjunction \emph{red or blue.} Section \ref{sssec:prod} likewise computed the product of copresheaves representing the conjunction \emph{red and blue}, and Section \ref{sssec:implication} computed the copresheaf representing the implication \emph{red} $\Rightarrow$ \emph{blue}.  In this section, we discuss the analogous constructions for \emph{enriched} copresheaves.  We have asserted several times that the $[0,1]$-category of $[0,1]$-copresheaves on a $[0,1]$-category $\mathcal{C}$ contains all of the enriched versions of (co)limits.  Without having yet said what the enriched version of (co)limits  are, one can understand that the reason $[0,1]$-enriched copresheaf categories are (co)complete is that the appropriate (co)limits are computed ``pointwise'' in $[0,1]$.  The real-analysis completeness of the interval implies that the infimum and supremum of any subset of the interval exists, hence the $[0,1]$-enriched categorical (co)limit of any diagram in $[0,1]$ exists.  This is analagous to the way (co)limits in the functor category $\mathsf{Set}^{\mathsf{C}}$ when $\mathsf{C}$ is an ordinary category are built up from (co)limits of sets, which always exist in $\mathsf{Set}$.

To make sense of (co)limits in this new enriched setting, the familiar definition must be slightly modified. This leads to the notion of \emph{weighted limits and colimits}, which are the appropriate notions of limits and colimits in enriched category theory.

\subsection{Weighted limits and colimits}\label{app:weighted_limits}
The appropriate notion of limits and colimits in enriched category theory are called \emph{weighted limits and colimits}. We focus on building intuition around limits as the discussion for colimits is analogous. Now,  to understand weighted limits, it will first help to recall the defining isomorphism for ordinary limits, also called ``conical limits.'' To that end, let $\mathsf{J}$ be an indexing category and let $\mathsf{C}$ be any category. Recall that the \emph{limit} of a diagram $F\colon \mathsf{J}\to\mathsf{C}$, if it exists, is an object $\lim F$ in $\mathsf{C}$ together with the following isomorphism of functors $\mathsf{C}\to\mathsf{Set}$,
\begin{equation}\label{eq:old_UP}
    \mathsf{C}(-,\lim F)\cong \mathsf{Set}^{\mathsf{J}}(\ast,\mathsf{C}(-,F))
\end{equation}
which is natural in the first variable. Here, $\ast$ is used to denote the constant functor at the one-point set:
\begin{equation}\label{eq:constant}
    \begin{tikzcd}
    \mathsf{J} \arrow[r,"\ast"] & \mathsf{Set}\\[-20pt]
    i \arrow[r, maps to]\arrow[d] & \ast\arrow[d,"\text{id}_\ast"]\\
    j \arrow[r, maps to] & \ast
    \end{tikzcd}
\end{equation}
So for each object $Z$ in $\mathsf{C}$ there is a bijection $\mathsf{C}(Z,\lim F)\cong \mathsf{Set}^{\mathsf{J}}(\ast,\mathsf{C}(Z,F))$
and so ``morphisms into the limit are the same as a cone over the diagram, whose legs commute with morphisms in the diagram,'' and one has in mind the following picture.
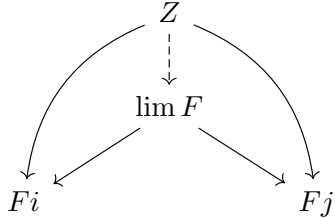
\[
\begin{tikzcd}[ampersand replacement=\&]
    \& Z \arrow[d, dashed] \arrow[ldd, bend right] \arrow[rdd, bend left] \& \\ \& \lim F  \arrow[dl] \arrow[dr]\&  \\ Fi \& \& Fj
\end{tikzcd}
\]
Taking a closer look at the right-hand side of the isomorphism in (\ref{eq:old_UP}), notice that for each object $Z$ in $\mathsf{C},$ a natural transformation from $\ast$ to $\mathsf{C}(Z,F)$ consists of a function $\ast=\ast(i)\to\mathsf{C}(Z,Fi)$ for each object $i$ in $\mathsf{J}$. This simply picks out a morphism $Z\to Fi$ in $\mathsf{C}$, and these morphisms comprise the legs of the limit cone over $\lim F$. Naturality ensures that these legs are compatible with the morphisms in the diagram.

With this background in mind, we now introduce weighted limits. The essential difference between the two constructions is that the constant functor in (\ref{eq:constant}) is replaced with a so-called \emph{$\mathcal{V}$-functor of weights} $W\colon\mathcal{J}\to\mathcal{V}$, where $\mathcal{V}$ is the base category over which enrichment takes place and where $\mathcal{J}$ is an indexing $\mathcal{V}$-category.  Here is the formal definition.

\begin{definition}
Let $\mathcal{V}$ be a closed commutative monoidal category and let $\mathcal{J}$ and $\mathcal{E}$ be $\mathcal{V}$-categories. Given a $\mathcal{V}$-functor $F\colon\mathcal{J}\to\mathcal{E}$ and a $\mathcal{V}$-functor $W\colon\mathcal{J}\to\mathcal{V}$, the \emph{weighted limit} of $F$, if it exists, is an object $\lim_WF$ of $\mathcal{E}$ together with the following isomorphism of $\mathcal{V}$-functors $\mathcal{E}\to\mathcal{V}$,
\begin{equation}\label{eq:new_UP}
    \mathcal{E}(-,\lim_WF)\cong\mathcal{V}^\mathcal{J}(W,\mathcal{E}(-,F))
\end{equation}
that is natural in the first variable.
\end{definition}
So in other words, for each object $Z$ in $\mathcal{E}$ there is an isomorphism $\mathcal{E}(Z,\text{lim}_WF)\cong\mathcal{V}^\mathcal{J}(W,\mathcal{E}(Z,F))$ of objects in $\mathcal{V}$.  The idea behind a  weighted colimit is analogous.  Given a $\mathcal{V}$-functor $F\colon\mathcal{J}\to\mathcal{E}$ and a $\mathcal{V}$-functor of weights $W\colon\mathcal{J}^{\text{op}}\to\mathcal{V}$, 
the \emph{weighted colimit} of $F$ is an object $\colim_WF$ of $\mathcal{E}$ together with an isomorphism 
\begin{equation}\label{eq:new_UP_colim}
    \mathcal{E}(\colim_WF,-)\cong\mathcal{V}^{\mathcal{J}^{\text{op}}}(W,\mathcal{E}(F,-))
\end{equation}
Details may be found in \cite[Chapter 7]{riehl_2014}.  

\subsection{Weighted products in $\widehat{\mathcal{L}}$}\label{ssec:prods}
Let's unwind the isomorphism in (\ref{eq:new_UP}) in the simple case when $\mathcal{V}=[0,1]$, when $\mathcal{E}=\widehat{\mathcal{L}}:=[0,1]^{\mathcal{L}}$, and when the indexing category $\mathcal{J}$ is a discrete category with two objects, call them $1$ and $2$, enriched over $[0,1]$ by setting $\mathcal{J}(i,j)=\delta_{ij}$ for $i,j\in\{1,2\}.$  To begin, fix a functor of weights $W\colon\mathcal{J}\to[0,1]$. This is nothing more than a choice of two numbers $w_1:=W(1)$ and $w_2:=W(2)$.  Further, for a fixed pair of copresheaves $f,g:\mathcal{L}\to [0,1]$, define $F\colon \mathcal{J}\to \widehat{\mathcal{L}}$ to be the $[0,1]$-functor with  $f:=F(1)$ and $g:=F(2)$.  
\begin{definition}Denote the weighted limit of $F$ with respect to the weight $W$ by \[(w_1 ,f)\times (w_2, g) := \lim_W F.\]
\end{definition}
\begin{theorem}The weighted limit
$(w_1, f)\times (w_2, g) :\mathcal{L} \to [0,1]$ is given by 
$\displaystyle
c \mapsto \min\left \{\frac{fc}{w_1},\frac{gc}{w_2},1 \right\}.$
\end{theorem}
\begin{proof}
To check that $c \mapsto \min\left \{\frac{fc}{w_1},\frac{gc}{w_2},1 \right\}$ satisfies the universal property of the weighted limit, let $Z\colon \mathcal{L} \to [0,1]$ be any copresheaf and look at Equation \eqref{eq:new_UP} evaluated at $Z$.  We need to check that
\begin{equation}\label{eq:new_UP_prod}
    \widehat{\mathcal{L}}(Z,\text{lim}_WF)\cong[0,1]^\mathcal{J}(W,\widehat{\mathcal{L}}(Z,F)).
\end{equation}
On the left hand side, with the claimed copresheaf substituted for the limit, we have
\[
\inf_{ c \in \mathcal{L}} \left\{\left [Zc,  \min\left \{\frac{fc}{w_1},\frac{gc}{w_2},1 \right\}\right]\right\} = \inf_{c \in \mathcal{L}}\left\{\frac{fc}{w_1Zc},\frac{gc}{w_2Zc},1\right\}.\]
Let's begin looking at the right hand side of Equation \eqref{eq:new_UP_prod} by simplifying the expression $\widehat{\mathcal{L}}(Z,F)$, which is a functor $J \to [0,1]$.  Evaluating at $i$ in $J$ yields the number 
\[\widehat{\mathcal{L}}(Z,Fi)=\inf_{c \in \mathcal{L}}\{[Zc,Fic]\} = \inf_{c \in \mathcal{L}}\left\{\frac{Fic}{Zc},1\right\}.\]
Now, the right hand side of Equation \eqref{eq:new_UP_prod}, as a hom-object in $[0,1]^\mathcal{J}$, is the minimum over the objects of $\mathcal{J}$ which are $1$ and $2$.  Using $f=F(1)$ and $g=F(2)$, we have:
\begin{align*}
[0,1]^\mathcal{J}(W,\widehat{\mathcal{L}}(Z,F))&=
\min\left\{\left[ w_1,\inf_{c \in \mathcal{L}}\left\{\frac{fc}{Zc},1\right\}\right],\left[ w_2,\inf_{c \in \mathcal{L}}\left\{\frac{gc}{Zc},1\right\}\right]\right\}\\
&=
\inf_{c \in \mathcal{L}}\left\{\frac{fc}{w_1Zc},\frac{gc}{w_2Zc},1\right\}.
\end{align*}
It remains to check that the assignment $c \mapsto \min\left \{\frac{fc}{w_1},\frac{gc}{w_2},1 \right\}$ is indeed a $[0,1]$-copresheaf; that is, that $\mathcal{L}(c,d)\leq \left[\lim_WF(c),\lim_WF(d)\right]$, or equivalently $\mathcal{L}(c,d)\lim_WF(c)\leq \lim_WF(d)$, for all $c,d\in\mathcal{L}$. This desired inequality follows from the simple observation that,
\begin{align*}
\mathcal{L}(c,d)\lim_WF(c) 
&= \min\left\{ \frac{\mathcal{L}(c,d)f(c)}{w_1},   \frac{\mathcal{L}(c,d)g(c)}{w_2},\mathcal{L}(c,d) \right\} \\
&\leq \min\left\{ \frac{f(d)}{w_1} , \frac{g(d)}{w_2},1\right\}\\
&\leq \min\left\{ \frac{f(d)}{w_1},\frac{g(d)}{w_2},1\right\}\\
&=\lim_WF(d),
\end{align*}
where the first inequality follows from the assumption that $f$ and $g$ are $[0,1]$-copresheaves.
\end{proof}
Keeping in mind that products in intuitionist logic serve as a kind of conjunction, let's look closer at the weighted product $(w_1, f)\times (w_2 , g)$ in the case $f$ and $g$ are representable.  Fix a pair of expressions $x$ and $y$ in $\mathcal{L}$ and a pair of nonzero weights $w_1, w_2 \in [0,1]$.  The weighted product $(w_1,h^x)\times (w_2,h^y)\colon\mathcal{L}\to[0,1]$ assigns to an expression $c$:
\begin{align}
((w_1, h^x)\times (w_2, h^y))(c)=\min\left\{\frac{h^x(c)}{w_1}, \frac{h^y(c)}{w_2},1\right\}.
\end{align} 
Remembering Equation (\ref{eq:enriched_representable}), the value of this minimum depends on whether the expressions $x$ and $y$ are jointly contained within $c$, and whether $\pi(c\vert x) \leq w_1$ and $\pi(c\vert y)\leq w_2$.
\begin{align}\label{eq:product}
((w_1,h^x)\times (w_2,h^y))(c) &=
	\begin{cases}
	\min\left\{\frac{\pi(c\vert x)}{w_1},\frac{\pi(c\vert y)}{w_2},1\right\} & \text{ if } x\leq c \text{ and } y\leq c\\[5pt]
	0 & \text{otherwise}.
	\end{cases}
\end{align}
The support of this copresheaf thus coincides with the support of logical ``and'' in Boolean logic.  As the weights decrease, the values of the quotients in Equation \eqref{eq:product} increase and thus contribute less to the value of $w_1h^x(c)\times w_2h^y(c)$, which is a minimum.  The copresheaf $(w_1,h^x)\times (w_2,h^y)$ captures something like a weighted conjunction.  

It's worth looking closer when the weights are $w_1=w_2=1$.  In this case, we denote the weighted product simply by $f\times g$.  This weighted product when the weights both equal $1$ works a lot like an ordinary product in the sense that morphisms into the product correspond precisely to products of morphisms.  Here, ``morphisms into the product" means a particular $[0,1]$-object, and ``products of morphisms" means the product of $[0,1]$-objects in $[0,1]$, which recall is a minimum.  That is,
\begin{lemma}\label{lem:enrichedproduct}
For any copresheaves $f,g,h\colon\mathcal{L}\to [0,1]$ we have 
\[\widehat{\mathcal{L}}(h,f\times g)=\widehat{\mathcal{L}}(h,f)\times\widehat{\mathcal{L}}(h,g).\]
\end{lemma}
\begin{proof} We compute:
\begin{align*}
\widehat{\mathcal{L}}(h,f\times g) &= \inf_{c\in \mathcal{L}}\{[hc,(f\times g)c]\} \\
&= \inf_{c\in \mathcal{L}}\{[hc, \min\{fc,gc\}] \}\\
&= \inf_{c\in \mathcal{L}}\left\{\frac{fc}{hc},\frac{gc}{hc},1\right\} \\
&= \min\left\{\inf_{c\in \mathcal{L}}\left\{ \frac{fc}{hc},1\right \}
,\inf_{c\in \mathcal{L}}\left\{\frac{gc}{hc},1\right\},1\right\} \\
&=\widehat{\mathcal{L}}(h,f)\times \widehat{\mathcal{L}}(h,g)
\end{align*}
\end{proof}
Lemma \ref{lem:enrichedproduct} will be helpful when we discuss enriched implications.  But before we do, let's discuss weighted coproducts.

\subsection{Weighted coproducts in $\widehat{\mathcal{L}}$}\label{ssec:coprods}
Now, let's look at a simple weighted coproduct in $\widehat{\mathcal{L}}:=[0,1]^\mathcal{L}$. Again, let the indexing category $\mathcal{J}$ be the same discrete category with two objects.  Fix a functor of weights $W\colon\mathcal{J}^{\text{op}}\to[0,1]$ setting $w_1:=W(1)$ and $w_2:=W(2)$ and a diagram $F\colon \mathcal{J}\to \widehat{\mathcal{L}}$ with  $f:=F(1)$ and $g:=F(2)$.  
\begin{definition}Denote the weighted colimit of $F$ with respect to the weight $W$ by \[(w_1, f)\sqcup (w_2, g) := \colim_W F.\]
\end{definition}
\begin{theorem} The weighted colimit
$(w_1, f) \sqcup  (w_2, g) \colon\mathcal{L} \to [0,1]$ is given by 
$\displaystyle
c \mapsto \max\left \{w_1fc,w_2gc\right\}.$
\end{theorem}
\begin{proof}Note that $\mathcal{J}=\mathcal{J}^{\text{op}}$ since there are no nonidentiy morphisms in $\mathcal{J}$.  Let $Z\colon \mathcal{L}\to [0,1]$ be any copresheaf.  We need to show 
\begin{equation}\label{eq:new_UP_coprod}
    \widehat{\mathcal{L}}(\colim_WF,Z)=[0,1]^\mathcal{J}(W,\widehat{\mathcal{L}}(F,Z)).
\end{equation}
Substituting the claimed colimit into the left hand side and evaluating at an object $i$ in $\mathcal{J}$ yields
\begin{align*}
\widehat{\mathcal{L}}(\colim_WFi,Z) &= \inf_{c \in \mathcal{L}}\{[\max\left \{w_1fc,w_2gc\right\},Zc]\}\\
&= \inf_{c \in \mathcal{L}}\left\{\frac{Zc}{\max\left \{w_1fc,w_2gc\right\}},1 \right\}\\
&= \inf_{c \in \mathcal{L}}\left\{\frac{Zc}{w_1fc},\frac{Zc}{w_2gc}\right\}.
\end{align*}
Evaluating $\widehat{\mathcal{L}}(F-,Z)$ at $i=1,2$ yields
\[\widehat{\mathcal{L}}(Fi,Z)=\inf_{c\in \mathcal{L}}\{[Fic,Zc]\} =\inf_{c\in \mathcal{L}} \left\{ \frac{Zc}{Fic}, 1 \right\}.
\]
Using $f=F(1)$ and $g=F(2)$ the right hand side of Equation \eqref{eq:new_UP_coprod} is
\begin{align*}
[0,1]^\mathcal{J}(W,\widehat{\mathcal{L}}(F,Z))&=\min\left\{\left[w_1,\inf_{c\in \mathcal{L}} \left\{ \frac{Zc}{fc}, 1 \right\}\right],\left[w_2,\inf_{c\in \mathcal{L}} \left\{ \frac{Zc}{gc}, 1 \right\}\right]\right\}\\
&=\inf_{c\in \mathcal{L}} \left\{ \frac{Zc}{w_1fc}, \frac{Zc}{w_2gc}, 1 \right\}.
\end{align*}
It remains to check that the assignment $c \mapsto  \max\left \{w_1fc,w_2gc\right\}$ is indeed a $[0,1]$-copresheaf; that is, that $\mathcal{L}(c,d)\leq \left[\colim_WF(c),\colim_WF(d)\right]$, or equivalently $\mathcal{L}(c,d)\colim_WF(c)\leq \colim_WF(d)$, for all $c,d\in\mathcal{L}$. This desired inequality follows from the simple observation that
\begin{align*}
\mathcal{L}(c,d)\colim_WF(c) 
&= \max\left\{\mathcal{L}(c,d)w_1f(c), \mathcal{L}(c,d)w_2g(c)\right\} \\
&\leq \max\left\{ \frac{f(d)}{w_1} , \frac{g(d)}{w_2}\right\}\\
&=\colim_WF(d),
\end{align*}
where the middle inequality follows from the assumption that $f$ and $g$ are $[0,1]$-copresheaves. This proves the copresheaf defined by $c \mapsto \max\left \{w_1fc,w_2gc\right\}$ satisfies the universal property to be the claimed weighted coproduct.
\end{proof}

The support of the copresheaf $(w_1, f) \sqcup (w_2, g)$ coincides with the support of logical ``or,'' and this matches the intuition that colimits capture something of disjunction.  Here, as the weights decrease, the corresponding cofactors contribute less to the weighted coproduct.

Now, let's now turn to a discussion of implication in the enriched setting. 

\subsection{Enriched implication}
Recall from Equation \eqref{eq:cartclosed} in Section \ref{sssec:implication} that the internal hom between a pair of ordinary copresheaves $G,H\colon\mathsf{C}\to\mathsf{Set}$ is the copresheaf $[G,H]$ defined on objects by $c\mapsto \widehat{\mathsf{C}}(h^c\times G,H)$, which is the set of natural transformations from $h^c\times G$ to $H$. 

Replacing the base category of sets with the unit interval gives the enriched version.
\begin{definition}
\label{def:enriched_implication}For any $f,g\colon\mathcal{L}\to [0,1]$ let
$[f,g]\colon\mathcal{L}\to [0,1]$ be defined by 
\begin{equation*}
[f,g](c):=\widehat{\mathcal{L}}(h^c\times f,g).
\end{equation*}
where $h^c\times f$ is the product defined above Lemma \ref{lem:enrichedproduct}.
\end{definition}
\begin{lemma}For any $f,g\colon\mathcal{L}\to [0,1]$ the function
$[f,g]\colon\mathcal{L}\to [0,1]$ is, in fact, a $[0,1]$-copresheaf.
\end{lemma}
\begin{proof}
We need to show that $\mathcal{L}(c,d) \leq [[f,g](c),[f,g](d)]$, which is equivalent to showing that $\mathcal{L}(c,d) [f,g](c) \leq [f,g](d)$.  
Start with the fact that $g$ is a copresheaf to get $\mathcal{L}(c,d) \leq [g(c),g(d)]$ which is equivalent to $\mathcal{L}(c,d) g(c) \leq g(d)$ which by the enriched Yoneda Lemma is equivalent to $\mathcal{L}(c,d) \widehat{\mathcal{L}}(h^c,g) \leq \widehat{\mathcal{L}}(h^d,g)$.  This is, in turn, equivalent to $\mathcal{L}(c,d) \left( \widehat{\mathcal{L}}(h^c,g)\times \widehat{\mathcal{L}}(f,g) \right)\leq \widehat{\mathcal{L}}(h^d,g)\times \widehat{\mathcal{L}}(f,g)$.  Using Lemma \ref{lem:enrichedproduct}, we can rewrite this inequality as $\mathcal{L}(c,d)  \widehat{\mathcal{L}}(h^c\times f,g)\leq \widehat{\mathcal{L}}(h^d\times f,g)$, which is equivalent to $\mathcal{L}(c,d)  \leq [\widehat{\mathcal{L}}(h^c\times f,g), \widehat{\mathcal{L}}(h^d\times f,g)]=[[f,g](c),[f,g](d)]$ as needed.
\end{proof}
We have an enriched functor $(f \times -)\colon\widehat{\mathcal{L}}\to\widehat{\mathcal{L}}$ and Definition \ref{def:enriched_implication} presents $[f,-]$ which is enriched right adjoint to $-\times f$.  The enriched Yoneda Lemma (Theorem \ref{thm:enrichedyoneda}) says $\widehat{\mathcal{L}}(h^c,[f,g])=[f,g](c)$ and by definition, we have $[f,g](c)=\widehat{\mathcal{L}}(h^c\times f,g).$  The equality $\widehat{\mathcal{L}}(h^c,[f,g])=\widehat{\mathcal{L}}(h^c\times f,g)$ for all representables $h^c$ implies equality for all copresheaves $h$:
\[
\widehat{\mathcal{L}}(h\times f,g) = \widehat{\mathcal{L}}(h,[f,g]).
\]
That is, $[f,g]$ serves as an internal-hom for $\widehat{\mathcal{L}}$ making it enriched \emph{Cartesian closed}.  Now, let's take a closer look at this internal hom $[f,g]$ when $f$ and $g$ are representable. 
\begin{theorem}For any
$x,y\in \mathcal{L}$, we have 
\begin{equation}\label{eq:representable_implication}
[h^x,h^y](c) = \inf_{d \in \mathcal{L}}\left\{\frac{\pi(d|y)}{\min\{\pi(d|c) , \pi(d|x)\}},1\right\}.
\end{equation}
\end{theorem}
\begin{proof}
We compute
\begin{align*}
[h^x,h^y](c) &= \widehat{\mathcal{L}}(h^c\times h^x,h^y)\\
&= \inf_{d \in \mathcal{L}} \left\{\left[(h^c \times h^x)(d),h^y(d)\right]\right\}\\
&= \inf_{d \in \mathcal{L}}\left\{\frac{\pi(d|y)}{\pi(d|c) \times \pi(d|x)},1\right\}.\\
&= \inf_{d \in \mathcal{L}}\left\{\frac{\pi(d|y)}{\min\{\pi(d|c), \pi(d|x)\}},1\right\}.
\end{align*}
For the last equality, remember that the categorical product in the unit interval is the minimum.
\end{proof}
 
Note that if there is a text $c$ that contains $x$ but does not contain $y$, then $[h^x,h^y](c)=0$ for the infimum in Equation \eqref{eq:representable_implication} is realized when $d=c$: the numerator $\pi(c|y)=0$ and the denominator $\pi(c|c) \times \pi(c|x)= \pi(c|x) \neq 0$.  So, if $[h^x,h^y](c)\neq 0$, then within the context $c$, either $c$ does not contain $x$ or $c$ contains both $x$ and $y$, capturing a quantitative sort of implication.

\begin{definition}
For any expressions $x,y\in \mathcal{L}$, define the implication $x\Rightarrow y$ to be the copresheaf $[h^x,h^y]\colon\mathcal{L}\to [0,1].$
\end{definition}


\section{A metric space interpretation}\label{sec:metric_space}
It is sometimes desirable to work in a category $\mathcal{M}$ that is a slight variant of the syntax category $\mathcal{L}$.  One can get from $\mathcal{L}$ to $\mathcal{M}$ by applying the negative logarithm to each hom object.  First notice the set of nonnegative extended reals $[0,\infty]$ together with addition of real numbers  is a commutative monoid with unit 0, where $a+\infty:=\infty$ and $\infty+a:=\infty $ for all $a$. If one further specifies a morphism from $a$ to $b$ whenever $b\leq a$, then  $[0,\infty]$ is a commutative monoidal preorder.  As a category $[0,\infty]$, like the unit interval, is also monoidal closed as well as complete and cocomplete:  the internal hom is given by \emph{truncated subtraction}, $[a,b]:=\max\{b-a,0\}$, and the limit of any diagram is \emph{supremum} of the numbers in the diagram, while the colimit is given by the \emph{infimum}. In particular, the categorical product and coproduct are respectively given by $a\times b = \max\{a,b\}$ and $a\sqcup b = \min\{a,b\}$ for all $a$ and $b$ in $[0,\infty]$.  The function $-\ln\colon [0,1]\to[0,\infty]$  defines an isomorphism of commutative monoidal preorders, the inverse of which $[0,\infty]\to [0,1]$ is the map $a\mapsto \exp(-a)$, and both maps are continuous and co-continuous isomorphisms of categories.

By applying $-\ln$ to morphisms of $\mathcal{L}$ we thus obtain a new category 
$\mathcal{M}$, enriched over the commutative monoidal preorder $[0,\infty]$, having the same objects as $\mathcal{L}$ and where the hom object between any pair of expressions $x$ and $y$ is given by $\mathcal{M}(x,y):=-\ln\mathcal{L}(x,y)$. It is then straightforward to check that $0\geq \mathcal{M}(x,x)$ for all expressions $x$ and moreover that $\mathcal{M}(x,y)+\mathcal{M}(y,z)\geq \mathcal{M}(x,z)$ for all expressions $x,y,z$ thus showing that $\mathcal{M}$ is indeed a $[0,\infty]$-category.  As was the case in $\mathcal{L}$, both inequalities in $\mathcal{M}$ are in fact equalities.  

\begin{figure}
    \centering
    \begin{tikzpicture}
	\node[tinydot] (R) at (0,0) {\tiny{red}};
	\node[tinydot] at (1,0) {};
	\node[tinydot] at (.4,.7) {};
	\node[tinydot] at (.5,.2) {};
	\node[tinydot] at (-.9,-.3) {\shortstack{\tiny{red}\\\tiny{ruby}}};
	\node[tinydot] at (-.2,-.8) {};
	\node[tinydot] at (-.4,-.4) {};
	\node[tinydot] at (-.1,1.2) {};
	\node[tinydot] at (-.5,1) {};
	\node[tinydot] at (-.2,.5) {};
	\node[tinydot] at (-.6,.3) {};
	\node[tinydot] at (.3,-.6) {};
	
	\node[tinydot] (A) at (-2,1.3) {};
	\node[tinydot]  at (-1,1) {};
	\node[tinydot]  at (-2,.6) {};
	\node[tinydot] (B) at (5,0.7) {\shortstack{\tiny{red}\\\tiny{idea}}};
	\node[tinydot] (C) at (2,1.5) {};
	\node[tinydot] (D) at (-2.5,-1.5) {};
	\node[tinydot] (E) at (-5,-0.5) {};
	\node[tinydot] (F) at (3.5,-1) {};
	\node[tinydot] (G) at (1,-1.5) {};
    \end{tikzpicture}
    \caption{When viewing language as a generalized metric space, a text such as \textit{red} is identified with its corresponding $[0,\infty]$-copresheaf $d_\mathcal{M}(\emph{red},-)$. Expressions that are unlikely continuations of \textit{red} are therefore far away, whereas expressions that are more likely to be continuations of \textit{red} are closer to it.}
    \label{fig:metricspace}
\end{figure}
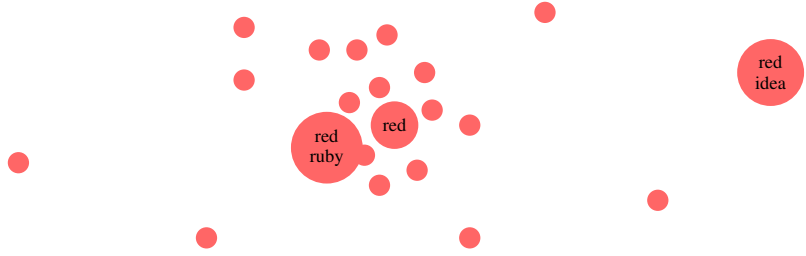

Categories enriched over $[0,\infty]$ are typically called \emph{generalized metric spaces} \cite{Lawvere73,Lawvere86} since composition of morphisms is precisely the triangle inequality, though notice that symmetry is not required as with  usual  metrics.  Even so, we embrace the generalized metric space point of view and will denote hom objects in $\mathcal{M}$ by $d_\mathcal{M}(x,y):=\mathcal{M}(x,y)$, thinking of $d_\mathcal{M}$ as defining a kind of distance between texts. Texts that are likely extensions of other texts have small distances from the texts they extend, and texts that are not extensions of one another are infinitely far apart. Figure \ref{fig:metricspace} illustrates a picture one might have in mind. Paths in this generalized metric space go in only one direction and represent stories:  you begin somewhere, there are expectations of where the story is going, the story continues, expectations of where you are going are revised, the story continues, and so on.  So  the two categories $\mathcal{L}$ and $\mathcal{M}$ have the same information, but the first emphasizes the probabilistic point of view while the second gives a geometric picture. 

Naturally, the idea that semantic information lies in a copresheaf category applies in $\mathcal{M}$ as well. There is a generalized metric space $\widehat{\mathcal{M}}:=[0,\infty]^\mathcal{M}$ whose objects are $[0,\infty]$-copresheaves on $\mathcal{M}$, which are functions $f\colon \mathcal{M}\to[0,\infty]$ satisfying $d_{M}(x,y)\geq [fx,fy]=\max\{fy-fx,0\}$ for all expressions $x$ and $y$ in the language.  We can think of $[0,\infty]$ as a generalized metric space itself and denote $[a,b]=\max\{b-a,0\}$ by $d_{[0,\infty]}(a,b)$.  So a copresheaf is like a metric contraction; that is, a function $f\colon\mathcal{M}\to [0,\infty]$ satsifying $d_{[0,\infty]}(fx,fy)\leq d_M(x,y)$.

Translating Equation (\ref{eq:hom_min}) tells us that the hom object between any pair of $[0,\infty]$-copresheaves $f$ and $g$ is given by 
\[\widehat{\mathcal{M}}(f,g)=\sup_{x\in\mathcal{M}}d_{[0,\infty]}(fx,gx)=\sup_{x\in\mathcal{M}}\max\{gx-fx,0\}\]
defining the generalized metric on the space of functions that one would expect.  Translating the results from Sections \ref{sec:enriched_copresheaves} and \ref{sec:enriched_prod_coprod}, one finds the generalized metric space $\widehat{\mathcal{M}}$ is also both complete and cocomplete with respect to weighted (co)limits, and moreover the Yoneda embedding   $\mathcal{M}^{\text{op}}\to\widehat{\mathcal{M}}$ that maps $x\mapsto d_\mathcal{M}(x,-)$ defines an isometric embedding of (the opposite category of) $\mathcal{M}$ as the representable copresheaves within $\widehat{\mathcal{M}}$.  Expressions $x$ and $y$ that may be unrelated in $\mathcal{M}$, can be combined in $\widehat{\mathcal{M}}$ in ways that are completely analagous to the weighted products and coproducts in $\widehat{\mathcal{L}}$.  See also \cite{willerton2013tight} for a similar discussion of generalized metric spaces and what is called the categorical Isbell completion.  

\begin{theorem}For any copresheaves $f,g\in \widehat{\mathcal{M}}$ and any weights $w_1, w_2$, the weighted product and coproducts of $f$ and $g$ are given by
\begin{gather*}
((w_1, f) \times (w_2, g))(c) = \max \{fc-w_1, gc-w_2,0\}  \\
((w_1, f) \sqcup (w_2, g))(c) = \min \{fc+w_1, gc+w_2\}
\end{gather*}
\end{theorem}
\begin{proof}
We translate the results from Sections \ref{ssec:prods} and \ref{ssec:coprods}.  Let $f,g\colon \mathcal{M} \to [0,\infty]$ and $w_1,w_2 \in [0,\infty]$.  Define $f',g'\colon\mathcal{L}\to [0,1]$ and weights $w_1',w_2'\in [0,1]$ by setting $f'c=\exp(-fc)$, $g'c=\exp(-gc)$, and $w_i'=\exp(-w_i)$ for $i=1,2$.  Then in $\widehat{\mathcal{L}}$ we have 
\[((w_1', f') \times (w_2', g') )(c)= \min\left\{\frac{f'c}{w'_1},\frac{g'c}{w'_2},1\right\}\]
Translating back to $\widehat{\mathcal{M}}$ by applying $-\ln$ yields
\begin{align*}
((w_1, f )\times (w_2, g))(c) &= -\ln\left( ((w_1', f') \times (w_2', g') )(c) \right)\\
&=-\ln \left(\min \left\{\frac{f'c}{w'_1},\frac{g'c}{w'_2},1\right\} \right)\\
&=\max \left\{-\ln\left(\frac{f'c}{w'_1}\right),-\ln\left(\frac{g'c}{w'_2}\right),0\right\}\\
&=\max \left\{-\ln(f'c)+\ln(w'_1),-\ln(g'c)+\ln(w'_2),0\right\}\\
&=\max \left\{fc-w_1,gc-w_2,0\right\}.\\
\end{align*}
Similarly,
\begin{align*}
((w_1, f) \sqcup (w_2, g))(c) &= -\ln\left( ((w_1', f') \sqcup  (w_2', g')) (c) \right)\\
&=-\ln \left(\max \left\{w'_1f'c,w'2g'c\right\} \right)\\
&=\min \left\{-\ln\left(w'_1f'c\right),-\ln\left(w'_2g'c\right)\right\}\\
&=\min \left\{-\ln(f'c)-\ln(w'_1),-\ln(g'c)-\ln(w'_2)\right\}\\
&=\min\left\{fc+w_1,gc+w_2\right\}.\\
\end{align*}
\end{proof}
Geometrically, one may think of nonrepresentable copresheaves, such as the weighted products and coproducts modeling conjoined and disjoined meanings, as new ``points'' added in  $\widehat{\mathcal{M}}$, which are specified precisely by giving their distance to all other points in the category $\mathcal{M}$.  The general idea that copresheaves $f$ in $\widehat{\mathcal{M}}$ are specified precisely by indicating the distance between $f$ and the representable copresheaves $d_\mathcal{M}(x,-)$ suggests that the categorical completion $\widehat{\mathcal{M}}$ resembles a kind of metric completion.  

\subsection{Tropical module structure}
We note now that the completion $\widehat{\mathcal{M}}$ has a semi-tropical module structure, that is, the structure of a module over the semi-tropical semi-ring that we define in the next paragraph.  The terminology \emph{semi}-tropical, and the connection to categorical cocompletions of generalized metric spaces can be found in \cite{willerton2013tight}.  In our context, the semi-tropical structure is analogous to the fact that the set of scalar valued functions on a given set has the structure of a vector space, which it inherits from the field of scalars. Here elements of $\widehat{\mathcal{M}}$ are functions valued in $[0,\infty]$, which is a semi-tropical semi-ring and therefore  $\widehat{\mathcal{M}}$ inherits the structure of a semi-tropical module.

\begin{definition}
The data $((-\infty,\infty],\oplus,\odot)$ with operations $\oplus$ and $\odot$ defined by \begin{equation*}
    s_1\oplus s_2=\min\{s_1,s_2\} \text{ and }s_1 \odot s_2=s_1+s_2
\end{equation*} 
defines a semi-ring called the \emph{tropical semi-ring} or the \emph{$(\min,+)$ algebra}.   
The sub semi-ring $([0,\infty],\oplus,\odot)$ is called the \emph{semi-tropical semi-ring}.  A \emph{module over a semi-ring} is a commutative monoid with an action of the semi-ring.
\end{definition}

\begin{theorem}The coproduct (with trivial weights $w_1=w_2=0$) in $\widehat{\mathcal{M}}$ makes it into a commutative monoid.  The map $[0,\infty]\times \widehat{\mathcal{M}} \to \widehat{\mathcal{M}}$ defined by $(s,f)\mapsto s\odot f$ where $(s\odot f)(x)=f(x)+s$ makes $\mathcal{M}$ into a module over the semi-tropical semi-ring.
\end{theorem}
\begin{proof}
We check: $(s_1 \oplus s_2)f (x) = \min\{s_1,s_2\}f(x)=\min\{fx+s_1,fx+s_2\}$ and $(s_1\odot s_2)f(x)=fx+s_1+s_2$ which is the same as $(s_1) \odot (s_2 \odot f )(x) = (s_2 \odot f )(x) +s_1 = fx+s_2+s_1.$
\end{proof}
In fact, the formulas for the weighted coproduct in $\widehat{\mathcal{M}}$ can be re-expressed as a tropical linear combination:
\begin{equation*}
((w_1, f) \sqcup ( w_2, g))(c) = \min\left\{fc+w_1,gc+w_2\right\} = w_1 \odot fc \oplus w_2 \odot gc.
\end{equation*}

We have a hunch that tropical geometry provides more than merely a different language to express the results described in $\widehat{\mathcal{L}}.$
Over the last twenty years, tropical geometry has been an area of active research.  One thing it provides is a way to transform non-linear algebraic geometry to piecewise linear geometry \cite{maclagan2012introduction, brugalle2015brief, Sturmfels, Sturmfels-erratum, yoshitomi2011generators, Sturmfels-Speyer1, Sturmfels-Speyer2, Viro}.  The passage is achieved by considering usual polynomials in several variables and 
replacing the additions and products by their tropical counterparts. The zero set of the original polynomial is mapped to the singular set of the tropical one which is a piecewise linear space.  \emph{Semi}-tropical module structures on presheaves on generalized metric spaces, along with a variation where $\min$ is replaced with $\max$, were studied in \cite{willerton2013tight}.

A potentially relevant discovery is a tropical module structure \cite{Sturmfels, Sturmfels-erratum,  Sturmfels-Speyer1} generated by distance functions, like the ones we have in this paper, but coming from a symmetric metric and satisfying
a tropical Pl\"ucker relation.  A one-dimensional tropical polytope arises, which realizes the metric defined by the distance functions as \emph{a tree metric}.  The achievement here is that extra limit points are added in precise locations, indicating branching points, so that the distances between leaves are realized as distances along tree paths connecting them.  This result is used in phylogenetics for the production of trees encoding common ancestors of species using distances defined by their DNA \cite{Sturmfels-Speyer1,Develin2004}.  We conjecture that a similar result for the metric semantic category $\widehat{\mathcal{M}}$ will give rise to phylogenetric structures that work like metric knowledge graphs.  

More speculatively, tropical structures might provide insights about how large language models actually learn semantic information.  Recently, it has been discovered  \cite{zhang2018tropical, Maragos1, Maragos2, Maragos3} that feedforward neural networks with ReLU activation functions compute tropical rational maps, opening a new way to study the mathematical structure underlying them.  LLMs use neural architectures to learn probabability distributions on text continuations.  That is, LLMs learn the syntax category $\mathcal{L}$.  It would be very interesting to determine if there is any link between the tropical module structures in the metric semantic language category $\widehat{\mathcal{M}}$ and some tropical rational maps computed within the neural architectures of LLMs.


\section{Conclusion}\label{sec:conclusion}
Knowing how to continue texts implies a good deal of semantic information.  So one might think that it would be very difficult to learn the syntactic language category $\mathcal{L}$ which encodes probability distributions on text continuations.  Two intelligent researchers debating five years ago could reasonably disagree about whether or not a computer model that coherently continues texts would require, as part of its training input, some semantic knowledge.  However, existing LLMs prove that is not the case.  Standard machine learning techniques demonstrate that the enriched category $\mathcal{L}$ can be learned in unsupervised way directly from samples of existing texts.  These models have then necessarily learned some semantic information.  In this paper, we provide a mathematical framework for where, mathematically speaking, that semantic information lives.  Namely, it lives in the category of $[0,1]$-copresheaves on $\mathcal{L}$.  We then described a number of categorical constructions that define meaningful operations on that semantic information.

\subsection{Applications and future directions}  There are several near-term applications of the work described in this paper.  One is to study the architectures of trained LLMs, like GPT3, which have successfully learned probability distributions on text continuations.  In the language of this paper, these distributions are precisely the representable enriched copresheaves.  One can look for other enriched categorical structures in these architectures and since tropical structures appear in both feedforward ReLU networks and the category of copresheaves, they might be a place to start.  Then the work in this paper could lead to well-defined mathematical operations on the parameters of trained language models that allow users to access, combine, and manipulate the semantic knowledge that is mathematically implicit in models that have learned to continue texts.  For a simple example, the concept of a gender neutral pronoun could be created by taking the coproduct of the copresheaves representing ``he'' and ``she.''  Also, weighted limits and colimits allow one to build concepts from multiple texts with precise shapes and weights, and by going back from copresheaves to texts, one would have control of the concepts contained in generated texts.  Realizing the internal hom operation, which is like a context-sensitive implication operator, would yield a powerful entailment tool, permitting certain texts to be input as given ``true'' when performing other NLP tasks.  Perhaps it is possible to use internal hom or the polyhedral structures arising from the tropical modules (like higher dimensional phylogenetic trees) to automatically create knowledge graphs, or otherwise organize semantic information, from the parameters of trained LLMs.  All of this depends on being able to realize certain categorical operations within existing architectures.  Another direction would be to design novel architectures which have the built-in capability to implement these categorical operations.  Using density operators to model probability distributions on text continuations is one idea \cite{BV2020}.  There are natural operations on density operators, like taking convex combinations, that correspond to operations on representable copresheaves.  One can average the density operators representing ``he'' and ``she'' to obtain a density that doesn't represent any particular word in the language, but rather captures the concept ``he'' $\sqcup$ ``she.''  Density operators also have spectral structures that can be accessed and manipulated realizing other operations with categorical interpretations.  Furthermore, tensor networks give highly efficient algorithms for storing and manipulating densities on existing classical hardware, making a tensor network language model an attractive possibility \cite[and references within]{DBLP:conf/aistats/MillerRT21}. 

In short, we have presented a mathematical framework that puts the kind of syntactical information that large language models learn into an enriched-category theoretic setting. Understanding that semantic information resides in a category of $[0,1]$-copresheaves and understanding how categorical operations act on that information leads to a number of concrete and appealing applications which should be explored further.

\section*{Acknowledgements}
The authors thank Olivia Caramello, Shawn Henry, Maxim Kontsevich, Laurent Lafforgue, Jacob Miller, David Jaz Myers, David Spivak, and Simon Willerton for helpful mathematical discussions.  The authors thank Juan Gastaldi and Luc Pellissier for discussions about their philosophical work \cite{Gastaldi1,doi:10.1080/03080188.2021.1890484} and the anonymous referees who made suggestions that greatly improved this article.

\nocite{Manning30046}
\bibliographystyle{alpha}
\bibliography{references}{}
\end{document}